\documentclass[a4paper]{article}  
     
\usepackage{graphicx}
\usepackage{amssymb,amsmath,mathrsfs,amsthm}
\usepackage{hyperref}
\usepackage{booktabs}
\usepackage{multirow}
\usepackage{algorithm}
\usepackage{algorithmic}
\usepackage{bm} % bold greek letters in math
\usepackage{breakurl}

\theoremstyle{plain}
\newtheorem{theorem}{Theorem}
\newtheorem{heuristic}{Heuristic}
\theoremstyle{definition}
\newtheorem{problem}{Problem}
\newtheorem{experiment}{Experiment}
\newtheorem{remark}{Remark}

\DeclareMathOperator{\End}{\mathrm{End}}

\DeclareMathOperator{\Cl}{\mathcal{CL}}

\DeclareMathOperator{\Exp}{\mathrm{E}}
\DeclareMathOperator{\Var}{\mathrm{Var}}
\DeclareMathOperator{\Cov}{\mathrm{Cov}}
\DeclareMathOperator{\Stdev}{\mathrm{Stdev}}

\newcommand{\F}{\mathbb{F}}
\newcommand{\Fbar}{\overline{\mathbb{F}}}
\newcommand{\OO}{\mathcal{O}}
\newcommand{\fa}{\mathfrak{a}}
\newcommand{\fb}{\mathfrak{b}}

\newcommand{\fl}{\mathfrak{l}}
\newcommand{\WB}{H}
\newcommand{\Hash}{H}   % Lower case $h$ already used for random walk steps.

\begin{document}

\title{Improved Algorithm for the Isogeny Problem \\ for Ordinary Elliptic Curves}

\author{Steven Galbraith \and Anton Stolbunov}

\date{May 31, 2011}

\maketitle

%%%%%%%%%%%%%%%%%%%%%%%%%%%%%%%%%%%%%%%%%%%%%%%%%%%%%%%%%%%%%%%%%%%%%%%%%%%%%%%
\begin{abstract}
A low storage algorithm for constructing isogenies between ordinary elliptic curves was proposed by Galbraith, Hess and Smart (GHS). We give an improvement of this algorithm by modifying the pseudorandom walk so that lower-degree isogenies are used more frequently. This is motivated by the fact that high degree isogenies are slower to compute than low degree ones. We  analyse the running time of the parallel collision search algorithm when the partitioning is uneven. We also give experimental results.
We conclude that our algorithm is around $14$ times faster than the GHS algorithm when constructing horizontal isogenies between random isogenous elliptic curves over a $160$-bit prime field.

The results apply to generic adding walks and the more general group action inverse problem; a speed-up is obtained whenever the cost of computing edges in the graph varies significantly. 
\end{abstract}

{\let\thefootnote\relax\footnotetext{E-mail addresses: \href{mailto:S.Galbraith@math.auckland.ac.nz}{S.Galbraith@math.auckland.ac.nz}, \ \href{mailto:anton@item.ntnu.no}{anton@item.ntnu.no}.}}

%%%%%%%%%%%%%%%%%%%%%%%%%%%%%%%%%%%%%%%%%%%%%%%%%%%%%%%%%%%%%%%%%%%%%%%%%%%%%%%
%%%%%%%%%%%%%%%%%%%%%%%%%%%%%%%%%%%%%%%%%%%%%%%%%%%%%%%%%%%%%%%%%%%%%%%%%%%%%%%
\section{Introduction}

Let $E_1$ and $E_2$ be elliptic curves over a finite field $\F_q$.
If $\#E_1( \F_q )  = \#E_2( \F_q )$
then there is an isogeny $\phi : E_1 \to E_2$ over $\F_q$ \cite[Theorem~1]{MR0206004}.
The \emph{isogeny problem} is to compute such an isogeny.
\begin{problem}[Isogeny Problem] 
Let $E_1/\F_q$ and $E_2/\F_q$ be ordinary elliptic curves satisfying $\#E_1( \F_q )  = \#E_2( \F_q )$. Compute an $\F_q$-isogeny $\phi : E_1 \to E_2$.
\end{problem}
The isogeny problem for ordinary elliptic curves (we do not consider the supersingular
case in this paper, though it is also interesting) over
finite fields is a natural problem, which has at least two important applications in cryptography. 

First, it allows to understand whether the difficulty of the discrete logarithm problem (DLP) is equal for all elliptic curves with the same number of points over $\mathbb{F}_q$. 
If $E_1$ and $E_2$ are ordinary then $\OO_1 = \End_{\Fbar_q}(E_1)$ 
and $\OO_2 = \End_{\Fbar_q}(E_2)$ are orders in a quadratic imaginary field $K$.
Let $\OO_K$ be the ring of integers of $K$ and define the \emph{conductor}
$c(E_i) = [ \OO_K : \OO_i ]$ for $i = 1,2$.
If there is a large prime $\ell$ such that $\ell \mid c(E_1)$ and $\ell \nmid c(E_{2})$
(or vice versa) then it seems to require at least $O( \ell^2 )$ steps to compute an
isogeny between $E_1$ and $E_2$, as explained in Section~\ref{sec:cost_isogenies}.
However, if this does not happen 
(in which case we say that the curves have \emph{comparable conductors})
then it can be feasible to compute an isogeny
from $E_1$ to $E_2$ using the algorithms due to  Galbraith~\cite{MR1728955} or  Galbraith, Hess and Smart~\cite{MR1975526} (GHS); the heuristic complexity is $\tilde{O}( q^{1/4 + o(1)} )$ bit operations.
As  has been observed by Jao, Miller and Venkatesan \cite{MR2236725}, and further discussed by Koblitz, Koblitz and Menezes \cite[\S11]{cryptoeprint:2008:390}, it follows
that the DLP  is random self-reducible among curves with the same number
of points and comparable conductors.

Second, the problem of constructing isogenies between ordinary elliptic curves is the basis of security of some recently proposed cryptographic schemes \cite{MR2210901,cryptoeprint:2006:145,cryptoeprint:2006:291,MR2654134,HeChenHu2011}. Cryptographic key sizes for these schemes should be chosen based on the complexity of the isogeny problem.

Galbraith, Hess and Smart~\cite{MR1975526} gave an algorithm, based on pseudorandom walks in the isogeny graph, to solve the problem.
At each step in the GHS algorithm an isogeny of relatively small degree $\ell$ is computed.
The starting point of our work is the observation that the cost of computing an isogeny
depends on~$\ell$ (see Fig.~\ref{fig:prime_action_timing}), and so it makes sense to choose a pseudorandom walk which ``prefers'' to use the fastest possible isogenies.
Similar ideas have also been used previously by authors:
Bisson and Sutherland~\cite{BissonSutherland09} in their algorithm for computing the endomorphism ring of ordinary elliptic curves; 
Stolbunov~\cite{MR2654134} in a family of cryptographic schemes based on isogenies.

The main problem is that making the pseudorandom walks ``uneven'' means that
the walks are ``less random'', and so the number of steps in the algorithm to solve
the isogeny problem increases.  However, this increase in cost is offset by the saving
in the cost of computing isogenies.
We analyse the effect of ``uneven'' partitions and suggest some good choices of parameters for the algorithm.
We also give experimental results to support our analysis.

The paper is organised as follows.
In Section~\ref{sec:2} we introduce a generalisation of the isogeny problem called the
\emph{group action inverse problem} (GAIP). We then explain why the isogeny problem
is the same as GAIP in the case of an ideal class group; we call this the
$\Cl$-GAIP.
In Section~\ref{sec:alg} we re-formulate (a variant of) the GHS algorithm as a generic algorithm for solving the GAIP and describe how it applies to the $\Cl$-GAIP. In Section~\ref{sec:RunningTime} we provide a theoretical analysis of the expected running time of the idealised algorithm.  Section~\ref{sec:prac} discusses how the idealised algorithm and the real implementation differ, and gives some experimental results. 
Section~\ref{sec:predictions} then makes some predictions about how the algorithm
will perform for isogeny computations, and determines the speedup of our ideas
compared with the algorithm described by Galbraith, Hess and Smart.
The main consequence of our work is that the isogeny problem can be solved
in less than one tenth of the time of the GHS algorithm.

%%%%%%%%%%%%%%%%%%%%%%%%%%%%%%%%%%%%%%%%%%%%%%%%%%%%%%%%%%%%%%%%%%%%%%%%%%%%%%%
%%%%%%%%%%%%%%%%%%%%%%%%%%%%%%%%%%%%%%%%%%%%%%%%%%%%%%%%%%%%%%%%%%%%%%%%%%%%%%%
\section{Definitions and Notation}\label{sec:2}
%%%%%%%%%%%%%%%%%%%%%%%%%%%%%%%%%%%%%%%%%%%%%%%%%%%%%%%%%%%%%%%%%%%%%%%%%%%%%%%
\subsection{The Group Action Inverse Problem}

Let $G$ be a finite abelian group, and $X$ a non-empty set. A (left) action of $G$ on $X$ is a map
\begin{align*}
G \times X &\to X\\
(g, x) &\mapsto g*x ,
\end{align*}
which satisfies the associativity property $(gh)*x = g*(h*x)$ for all $g, h \in G$, $x\in X$, and the property $e*x=x$ for the identity element $e\in G$ and all $x\in X$. The \emph{orbit} of a set element $x\in X$ is the subset $G*x=\{g*x \mid g\in G\}$. The orbits of the elements of $X$ are equivalence classes. The \emph{stabilizer} of $x$ is the set of all elements in $G$ that fix $x$: $G_x = \{g \in G \mid g*x = x\}$. 

\begin{problem}[Group Action Inverse Problem] 
Let $G$ be a finite abelian group acting on a non-empty set $X$. Given elements $x, y\in X$, find a group element $g\in G$ such that $g*x=y$.
\end{problem}

When the action of $G$ on $X$ is \emph{transitive}, that is, $X$ is finite and there is only one orbit, then the GAIP has at least one solution. When the action is \emph{free}, i.e. the stabilizer of any set element is trivial, then the GAIP has at most one solution. In the case of a free and transitive action, the set $X$ is called a \emph{principal homogeneous space} for the group $G$, and the GAIP has exactly one solution. This last type of GAIP will be considered in the rest of the paper.

%%%%%%%%%%%%%%%%%%%%%%%%%%%%%%%%%%%%%%%%%%%%%%%%%%%%%%%%%%%%%%%%%%%%%%%%%%%%%%%
\subsection{The Isogeny Problem and the Class Group Action Inverse Problem}
\label{sec:Cl-GAIP-and-isogeny}

Recall from the introduction that
$E_1$ and $E_2$ are ordinary elliptic curves over $\F_q$ with $\#E_1( \F_q )  = \#E_2( \F_q )$,
 $\OO_i = \End_{\Fbar_q}(E_i)$ and
$c(E_i) = [ \OO_K : \OO_i ]$ for $i = 1,2$.
As noted by Galbraith~\cite{MR1728955} (building on work of Kohel~\cite{Koh96}), a natural approach to compute an isogeny from
$E_1$ to $E_2$ is to first take ``vertical'' isogenies to elliptic curves $E_1'$ and $E_2'$
such that $\End_{\Fbar_q}(E_i') = \OO_K$, and the isogeny problem is reduced to computing a
``horizontal'' isogeny from $E_1'$ to $E_2'$.
Alternatively, if $\OO_1$ and $\OO_2$ are comparable, but both $c(E_1)$
and $c(E_2)$ have a large prime factor, one can use horizontal and/or vertical
isogenies from $E_1$ to a curve $E_1'$ such that $\End_{\Fbar_q}(E_1') = \OO_2$
and the problem is again reduced to computing a horizontal isogeny.

So, without loss of generality, we  assume for the remainder of the paper
that $\End_{\Fbar_q}(E_1) = \End_{\Fbar_q}(E_2) $.
Define $ \OO$ to be the order $\End_{\Fbar_q}(E_1)$.
Write $\Cl(\OO )$ for the group of invertible $\OO$-ideals modulo 
principal $\OO$-ideals and $h( \OO )$ for the order of $\Cl(\OO)$.

The theory of complex multiplication implies that there are $h(\OO)$ isomorphism
classes of elliptic curves $E$ over $\F_q$ with $\End_{\Fbar_q}(E) = \OO$ and a fixed number of points $\#E(\F_q)$.
There is a (non-canonical) one-to-one-correspondence between
isomorphism classes of elliptic curves $E $ over $\F_q$ with $\End_{\Fbar_q}(E) = \OO$
and ideal classes in $\Cl(\OO )$ \cite{MR0265369}.
There is a (canonical) one-to-one correspondence between invertible $\OO$-ideals
$\fl$ and isogenies, such that if $\fl$ is an ideal of norm $\ell$ 
and $E$ is an elliptic curve corresponding to the ideal $\fa$ then
there is an $\ell$-isogeny from $E$ to $E'$ where $E'$ corresponds to
the ideal $\fa \fl^{-1}$.
Galbraith, Hess and Smart~\cite{MR1975526} show how, given an elliptic curve $E$ and an ideal $\fb$,
one can efficiently compute an explicit isogeny $\phi : E \to E'$ corresponding to $\fb$ via
the above correspondence.

Let $X$ be the set of isomorphism classes of elliptic curves over $\F_q$ with $\End_{\Fbar_q}(E) = \OO$ and a fixed $\#E(\F_q)$. It follows that $\Cl( \OO )$ acts on $X$ and so we can define $\fb * E$ to be
the isomorphism class of the image curve for the isogeny corresponding to $\fb$.
The horizontal isogeny problem is a special case of the GAIP, which we call the class group action inverse problem ($\Cl$-GAIP).

\begin{problem}[Class Group Action Inverse Problem] 
Let $E_1/\F_q$ and $E_2/\F_q$ be ordinary elliptic curves satisfying $\#E_1( \F_q )  = \#E_2( \F_q )$ and $\End_{\Fbar_q}(E_1) = \End_{\Fbar_q}(E_2) = \OO$. Find the ideal class $[\fb]\in\Cl(\OO)$ such that the curves $\fb*E_1$ and $E_2$ are isomorphic.
\end{problem}
Hence, for the rest of the paper we study the GAIP, keeping in mind
this specific application.

Let $\WB = \{ \fl_1, \dots, \fl_r \}$ be a set of distinct prime ideals.
We define the \emph{ideal class graph} to be the graph with vertex
set $\Cl(\OO )$ and, for each $\fl \in \WB$, an edge $( \fa , \fa \fl^{-1} )$
for all $\fa \in \Cl( \OO )$.
Similarly, we define the \emph{isogeny graph} to have vertex set
being isomorphism classes of elliptic curves with endomorphism ring $\OO$
and an edge between two isomorphism classes if there is an isogeny between
them corresponding to an ideal $\fl \in \WB$.

%%%%%%%%%%%%%%%%%%%%%%%%%%%%%%%%%%%%%%%%%%%%%%%%%%%%%%%%%%%%%%%%%%%%%%%%%%%%%%%
\subsection{Other Notation}
By $a \leftarrow b$ we denote the assignment of value $b$ to a variable~$a$. By $a \xleftarrow{R} G$ we mean that $a$ is sampled from the uniform distribution on the set of elements of~$G$. We write $\#S$ for the number of elements in~$S$. By $\log (n)$ we denote the binary logarithm of~$n$. All equalities of the form $f(x)=O(g(x))$ are one-way equalities that should be read as ``$f(x)$ is $O(g(x))$''.

%%%%%%%%%%%%%%%%%%%%%%%%%%%%%%%%%%%%%%%%%%%%%%%%%%%%%%%%%%%%%%%%%%%%%%%%%%%%%%%
%%%%%%%%%%%%%%%%%%%%%%%%%%%%%%%%%%%%%%%%%%%%%%%%%%%%%%%%%%%%%%%%%%%%%%%%%%%%%%%
\section{Algorithm for Solving the GAIP and the \\\texorpdfstring{$\Cl$}{CL}-GAIP}
\label{sec:alg}
%%%%%%%%%%%%%%%%%%%%%%%%%%%%%%%%%%%%%%%%%%%%%%%%%%%%%%%%%%%%%%%%%%%%%%%%%%%%%%%
\subsection{Previous Isogeny Problem Algorithms}
\label{sec:GHS}

The first algorithm for solving the isogeny problem (equivaently, the $\Cl$-GAIP) was proposed by Galbraith~\cite{MR1728955}. 
Let $E_1$ and $E_2$ be elliptic curves over $\F_q$ with $\End( E_i ) = \OO$
(alternatively, let $x$ an $y$ be $\OO$-ideal classes).
The idea was to construct two graphs of elliptic curves (subgraphs of the
isogeny graph), one rooted at
$E_1$ and the other at $E_2$ (equivalently, two subgraphs of the ideal class graph rooted
at $x$ and $y$ respectively). Edges in the graph correspond to small-degree ideals.
By the birthday paradox, when the graphs have total size approximately $\sqrt{\pi h(\OO) }$ one expects them to have a vertex in common, in which case we have a path of isogenies from $E_1$ to $E_2$.
Indeed, under the assumption that the graphs behave like random subgraphs from the point of view of their intersection, it is natural to conjecture that the algorithm halts when the total number of vertices visited is, on average, $\sqrt{\pi h(\OO)}$.
Note that this algorithm requires an exponential amount of time and memory.

The second, and previously the best, algorithm was due to Galbraith, Hess and Smart~\cite{MR1975526} (in particular the stage~1 of the algorithm described in that paper). The major improvement was to use pseudorandom walks and parallel collision search in the isogeny graph, rather than storing entire subgraphs.
We give a generic description of this method in the next section.
The advantage of the GHS method is that
it only requires a polynomial amount of memory, and can be easily parallelised or distributed.

Although this paper considers the classical computational model, we note that a subexponential-time quantum algorithm for the isogeny problem has been proposed by Childs, Jao and Soukharev~\cite{quantisog}.

%%%%%%%%%%%%%%%%%%%%%%%%%%%%%%%%%%%%%%%%%%%%%%%%%%%%%%%%%%%%%%%%%%%%%%%%%%%%%%%
\subsection{Generic Description of the GAIP Solving Algorithm}
\label{sec:Alg}
Let the GAIP $(x,y)$ be defined for a group $G$ acting on a set $X$, and let $r$ be a positive integer greater or equal than the rank of $G$. Choose a generating set $\WB =\{g_1, \dots, g_r\}\subset G$ and consider a graph $\Gamma$ with vertices the elements of $X$, and edges $(z,g_i*z)$, for all $1\le i\le r$. In the special case $G = \Cl( \OO )$, $X$~the set of isomorphism classes of elliptic curves with the endomorphism ring~$\OO$, and $\WB = \{ \fl_1, \dots, \fl_r \}$, we obtain the isogeny graph defined in Section~\ref{sec:Cl-GAIP-and-isogeny}. 

To solve the GAIP it suffices to find an (undirected) path in $\Gamma$ between $x$ and $y$. A natural way to do this is to use (pseudo)random walks in $\Gamma$, starting from $x$ and $y$. For instance one can use a random function $v\colon X\to\{1, \dots, r\} $ and the map
\begin{align*}
\psi\colon X &\to X\\
z &\mapsto g_{v(z)}*z .
\end{align*}
The following language will be used throughout the paper:
the function $v(z)$ is a \emph{partitioning function}, because it defines a partition $P$ on the set $X$. By an abuse of notation we will call parts in $P$ \emph{partitions}. Note that we do not require all partitions to be of the same size. \emph{Partitioning probabilities} $p_1, \dots, p_r$ are defined as
\[ p_i=\Pr\left[v(z)=i\mid z\xleftarrow{R} X\right] \qquad \text{for all } 1\le i\le r . \]
A \emph{walk} on $\Gamma$ is a sequence of nodes computed as 
\begin{equation*}
z_{j+1}=\psi(z_j) .
\end{equation*}
A \emph{hop} is one edge in the graph (i.e., one step of the walk). The set $\WB$ is called \emph{the supporting set} for walks on $\Gamma$. The above walk is a generalization of the \emph{adding walk} proposed by Teske for groups~\cite{MR1697652}.

One can apply the parallel collision search concept, as proposed by van Oorschot and Wiener~\cite{MR1664774}. To do this, define a subset $X_D$ of \emph{distinguished elements} in $X$, such that it is easy to verify that $z\in X_D$. Pseudorandom walks in $\Gamma$ are formed by taking a random initial vertex\footnote{The GHS algorithm~\cite{MR1975526} does not specify how to sample random vertices in the isogeny graph. We use an algorithm from Stolbunov~\cite[\S6.1]{MR2654134}, which will be briefly explained at the end of Section~\ref{sec:better-choices}.}, moving along edges with a certain probability, and halting when the current vertex is a distinguished element. This framework was used by Galbraith, Hess and Smart~\cite{MR1975526}. Figure~\ref{fig:A} presents Algorithm $\mathcal{A}$, which is an algorithm to solve the GAIP following this approach.

Algorithm $\mathcal{A}$  uses $2t$ client threads, where $t\ge 1$, and one server thread. The algorithm takes as input a GAIP instance $(x_0,x_1)$ and an integer $t$. The server starts $t$ clients, each performing a walk starting from a randomized node $h_{0,i}*x_0$ for $1 \le i \le t$.
The server starts another $t$ clients, each performing a walk starting from a randomized node $h_{1,i}*x_1$. Each client continues the deterministic pseudorandom walk until it hits a distinguished node. Once a thread hits a distinguished node $z=a*x_s$, it puts the triple $(z, a, s)$ on the shared queue  and terminates. The server stores all received triples in a database $D$ and restarts clients from new randomized starting nodes.

\begin{figure}
\begin{minipage}[t]{0.47\textwidth}
\begin{algorithm}[H]
\caption{Server} 
\label{alg:Aserver} 
\begin{algorithmic}[1]
\REQUIRE $(x_0,x_1,t)\in X\times X\times \mathbb{N}$
\FOR{$i=1$ to $t$}
\STATE $(h_0,h_1)\xleftarrow{R}G\times G$
\STATE start $\mathrm{client}(h_0*x_0,h_0,0)$
\STATE start $\mathrm{client}(h_1*x_1,h_1,1)$
\ENDFOR
\STATE $D\leftarrow\{\}$
\WHILE{true}
\STATE fetch $(z,a,s)$ from queue
\IF{$(z,b,1-s)\in D$ for some $b$}
\STATE break loop
\ENDIF
\STATE $D\leftarrow D \cup \{(z,a,s)\}$
\STATE $h\xleftarrow{R}G$
\STATE start $\mathrm{client}(h*x_s,h,s)$
\ENDWHILE
\STATE stop all clients
\ENSURE $a^{1-2s}b^{2s-1}$
\end{algorithmic}
\end{algorithm}
\end{minipage}
\hfill
\begin{minipage}[t]{0.47\textwidth}
\begin{algorithm}[H]
\caption{Client} 
\label{alg:Aclient} 
\begin{algorithmic}[1]
\REQUIRE $(z,a,s)\in X\times G\times \{0,1\}$
\STATE $c\leftarrow 0$
\WHILE{$z\not\in X_D$}
\STATE $i\leftarrow v(z)$
\STATE \label{state:l4} $z\leftarrow g_i * z$
\STATE $a\leftarrow a g_i$
\STATE $c\leftarrow c+1$
\IF{$c>c_{\max}$} \label{state:l7}
\STATE $(z,a,s)\leftarrow\perp$
\STATE break loop
\ENDIF \label{state:l10}
\ENDWHILE
\ENSURE $(z,a,s)$
\end{algorithmic}
\end{algorithm}
\end{minipage}
\caption{Algorithm $\mathcal{A}$ for solving the GAIP.}\label{fig:A}
\end{figure}

A \emph{collision} is an event when some node is visited by client threads twice, while the preceding nodes visited by the threads are different. Since the walks are deterministic, after a collision the two threads follow the same route unless they hit a distinguished node. Thus every collision results in two triples of the form $(z,\cdot,\cdot)$ being submitted to the server. A collision of walks, one of which was started from $x_0$ and the other one from $x_1$, is called a \emph{good collision}. After a good collision the server detects two triples $(z,a,0)$ and $(z,b,1)$. It then halts all clients and outputs the solution $b^{-1}a$. 

Since a walk might loop before it hits a distinguished node, clients use a simple loop detection mechanism that checks whether the walk remains shorter than a fixed maximum length $c_{\max}$. The value $c_{\max}$ is usually chosen to be a function of $\theta$, e.g. $c_{\max}=30/\theta$, which means that walks $30$ times longer than expected are abandoned\footnote{Van Oorschot and Wiener~\cite{MR1664774} suggest $c_{\max}=20/\theta$. Our value is larger in order to preserve more non-looped walks.}.

Denote by $\alpha$ the number of nodes visited by Algorithm $\mathcal{A}$, counted with repetition. If nodes were sampled uniformly at random then the expected value $\Exp(\alpha)$ would be close to $\sqrt{ \pi \#G }$ by a variant of the birthday paradox (see Section~\ref{sec:prevres}). The expected total (serial) running time of Algorithm $\mathcal{A}$ approximately equals the product of $\Exp(\alpha)$ with the average cost of computing $g_i * z$ in line~\ref{state:l4} of the client algorithm\footnote{We do not count database access times and expected $L \theta \sqrt{n}$ random samplings of a group element.}. Our main observation is that the cost of computing $g_i * z$ is not the same for all $g_i$. Hence, one can speed up the algorithm by favoring the $g_i$ which are faster to compute.

In the $\Cl$-GAIP, the supporting set $\WB$ is usually chosen to consist of prime ideals above the smallest integer primes which split in $\OO$. In some rare cases it may be necessary to add one or more prime ideals of larger norm to ensure that $\WB$ generates $\Cl( \OO )$. Ramified primes can also be used, but since their order equals two in $\Cl( \OO )$ they suffer from the defect mentioned in the next section.

%%%%%%%%%%%%%%%%%%%%%%%%%%%%%%%%%%%%%%%%%%%%%%%%%%%%%%%%%%%%%%%%%%%%%%%%%%%%%%%
\subsection{A Remark on the GHS Algorithm}

The GHS paper~\cite{MR1975526} states that ``it is usually enough that $\WB$ contains about $16$ distinct split primes'', and the partitioning function should ``have a distribution close to uniform''. In other words, it was advised to use about $r = 16$ partitions of approximately equal size. 
We will compare our algorithm against those suggested parameters in the remainder of the paper.

We note a potentially serious problem\footnote{This remark also applies to the isogeny walk given by Teske \cite[Algorithm~1]{MR2210901}. Interestingly, another isogeny walk is given in Algorithm~3 of the same paper, which is not affected by this problem.} with the algorithm of Galbraith, Hess and Smart~\cite{MR1975526}. On every hop the algorithm chooses a small prime $\ell$ and a bit $b$ uniformly at random. Typically,  $\ell$ is split and the algorithm chooses one of the two $\ell$-isogenous elliptic curves deterministically using the bit $b$. Hence for a fixed $\ell$, every hop where $\ell$ is chosen produces an action by, equally likely, the ideal $\mathfrak{l}$ or $\mathfrak{l}^{-1}$ (where
$(\ell)  = \mathfrak{l} \mathfrak{l}^{-1}$). Thus, since the ideal class
group is abelian, the expected power of the ideal $\mathfrak{l}$ that has acted on the starting elliptic curve after any number of hops equals $0$. Such a walk is far from random, as it tends to remain ``close'' to its initial node. 
Hence, most likely the method of Galbraith, Hess and Smart does not perform
as well in practice as the heuristic predictions stated in~\cite{MR1975526}.
To avoid this problem, our algorithm always acts by the same ideal $\mathfrak{l}$ when the prime $\ell$ is chosen (i.e., the set $\WB$ never contains both $\fl$ and $\fl^{-1}$; unless $\fl$ is ramified). We stress that the speed improvement of our algorithm is not due to the correction of the named flaw but because of the use of an uneven partitioning.

%%%%%%%%%%%%%%%%%%%%%%%%%%%%%%%%%%%%%%%%%%%%%%%%%%%%%%%%%%%%%%%%%%%%%%%%%%%%%%%
\subsection{Better Choices for Solving the \texorpdfstring{$\Cl$}{CL}-GAIP}
\label{sec:better-choices}

We now discuss the main idea of the paper, which is to make the pseudorandom walks faster by using smaller degree prime ideals more often than larger degree ones.

Recall that $\alpha$ denotes the number of nodes visited by Algorithm $\mathcal{A}$, counted with repetition, and that $\Exp(\alpha)$ is close to $\sqrt{ \pi n }$, where $n=\#G$. Therefore it is more convenient to consider the variable 
\[ L=\frac{\alpha}{\sqrt{n}} . \]
The value of $L$ is fully determined by 
the group, 
the problem instance $(x_0,x_1)$,
the supporting set $\WB$,
the partitioning function $v()$, 
the subset $X_D$ of distinguished nodes,
the loop detection value $c_{\max}$
and the random choices made by the algorithm.
We define $\Exp(\,L\mid r,\vec p, m , \theta, c_{\max}\,)$ to be the expected value of~$L$, taken over random choices of all the above parameters, conditioned on the values of the parameters:
\begin{description}
\setlength{\itemsep}{1pt}\setlength{\parskip}{1pt}
\item[$r$] the number of partitions;
\item[${\vec p =  (p_1,\dots,p_r)}$] the partitioning probabilities;
\item[$m = \lceil \log (n)\rceil $] the ceiling function of the binary logarithm of $\#G$;
\item[$ \theta $] the probability of distinguished nodes;
\item[$ c_{\max} $] the loop detection value.
\end{description}
To shorten the notation we will write $\Exp(L)$ instead of $\Exp(\,L\mid r,\vec p, m , \theta, c_{\max}\,)$.

The average running time of a step in the algorithm (equivalently, hop) is $\vec p\, \vec t = \sum_{i=1}^r p_i t_i$, where $\vec t$ is a column vector of timings of actions by the $r$ chosen primes (see Fig.~\ref{fig:prime_action_timing} for such timings).
Hence, the expected serial running time of Algorithm $\mathcal{A}$ is approximately 
\begin{equation}
\label{eq:exptime}
\Exp(L)\,\sqrt{n}\, \vec p\, \vec t . 
\end{equation}

Ideally, the number of partitions $r$ and the probability distribution $\vec p$ should be chosen by solving the optimization problem: given $n, \theta, \vec t$, choose $r$ and $\vec p$ to minimise the expected running time $\Exp(L)\,\sqrt{n}\,\vec p\, \vec t$. We do not claim in this paper a complete solution to this optimisation problem. But we do discuss how $\Exp(L)$ depends on $r$ and $\vec p$, and we suggest some choices for these parameters.

For simplicity, and because they seem to give good results in practice, we restrict our attention to vectors $\vec p = (p_1, \dots, p_r )$ such that the probabilities are in geometric progression $p_{i+1}/ p_{i} = w$ for $1 \le i < r$. For example, taking $r = 4$ and $w=1/2$ means probabilities $( p_1, \tfrac{1}{2} p_1, \tfrac{1}{4} p_1, \tfrac{1}{8} p_1 )$ which add up to $1$ (and so $p_1 = 8/15 \approx 0.53$). In our practical analysis we restrict to $3\le r\le 16$ and $\vec p$ is the geometric progression of ratio $w\in\{1,3/4,1/2,1/3,1/4\}$. This choice is probably not the best solution to the optimization problem, but it seems to work well in practice.

To implement the starting randomization of the walks we use a method proposed by Stolbunov~\cite[\S6.1]{MR2654134}. We briefly describe the method. Since the class group structure computation is much faster~\cite{MR2654130} than Algorithm~$\mathcal{A}$, one first computes the class group structure. For an imaginary quadratic order $\OO$ of discriminant $\Delta$, the class group $\Cl(\OO)$ is generated (assuming GRH) by the set $\mathcal{L}$ of prime ideals of split norms less than or equal to $\ell_{\max}=c_1 \log^2\lvert\Delta\rvert$, for an effectively computable constant $c_1$ \cite[Corollary~6.2]{MR702519}. Note that the set $\mathcal{L}$ used for the random sampling can be larger than the supporting set $\WB$. Knowing the class group generators and their orders, one obtains a random group element in a smooth form by raising generators to random exponents, each chosen between zero and the corresponding order. To shorten the representation one reduces it modulo the lattice of relations among the elements of~$\mathcal{L}$.
Indeed, it is possible to write any element of $\Cl(\OO)$ as an $O(\log\lvert\Delta\rvert)$-term product of elements in $\mathcal{L}$. Jao, Miller and Venkatesan have shown (assuming GRH) that the ideal class graph $(\Cl(\OO), \mathcal{L})$ is an expander graph~\cite[Theorem 1.5]{MR2521489}. Since the diameter of an expander graph is less than or equal to $2\log (h)/\log(1+c)$ for the expansion coefficient $c$ and the number of vertices $h$ \cite[Theorem~9.9]{Goldreich2001LN}, the diameter of the ideal class graph $(\Cl(\OO),\mathcal{L})$ is $O(\log (h))$, where $h\approx\lvert\Delta\rvert^{1/2}$. 

%%%%%%%%%%%%%%%%%%%%%%%%%%%%%%%%%%%%%%%%%%%%%%%%%%%%%%%%%%%%%%%%%%%%%%%%%%%%%%%
%%%%%%%%%%%%%%%%%%%%%%%%%%%%%%%%%%%%%%%%%%%%%%%%%%%%%%%%%%%%%%%%%%%%%%%%%%%%%%%
\section{Theoretical Analysis of the Algorithm}
\label{sec:RunningTime}
%%%%%%%%%%%%%%%%%%%%%%%%%%%%%%%%%%%%%%%%%%%%%%%%%%%%%%%%%%%%%%%%%%%%%%%%%%%%%%%
\subsection{Previous Results}
\label{sec:prevres}
A tremendous amount of research on the running time analysis of the Pollard rho algorithm has been carried out by various authors. We give a brief overview of some of the results relevant to our work.

First we consider random mappings on a set $X$ of $n$ elements. Rapoport~\cite[\S II]{springerlink:10.1007/BF02477489} and Harris~\cite[\S 3]{MR0119227} obtained an approximation for the expected value of the number $\rho$ of distinct elements in a random walk on~$X$: 
\[
\Exp(\rho)\approx\sqrt{\frac{\pi n}{2}} .
\]
For a more precise statement see Knuth~\cite[Exercise 3.1.12]{KnuthV2}.
These results were subsequently used to approximate the expected length of the rho-shaped walk in the Pollard's algorithm~\cite{MR0392798}.

Van Oorschot and Wiener~\cite[\S 4.1]{MR1664774} proposed a parallel version of the Pollard's rho algorithm. When more than one walk is run in parallel, several collisions can occur, and only some of them may be useful (we call these collisions good). Let $\mathfrak{p}$ be the probability that a random collision is good. They obtained the following approximation for the expected value of the number $\lambda$ of distinct visited nodes, when the number of collisions is small:
\begin{equation} \label{eq:vOW}
\Exp(\lambda)\approx\sqrt{\frac{\pi n}{2 \mathfrak{p}}} .
\end{equation}

The iteration function proposed by Pollard~\cite{MR0491431} for the DLP involved three partitions of approximately equal size: two corresponding to multiplication and one to squaring hops. Teske proposed a different type of iteration function which she called an adding walk~\cite{MR1697652}. Adding walks allowed more partitions, but it was still preferable to have equally-sized partitions because the costs of iterations were approximately equal. 
Brent and Pollard~\cite{MR606520} and
Blackburn and Murphy~\cite{BlaMur98} provided a heuristic argument where they assumed that the restrictions of the iterating function to $r$ equally-sized partitions were random mappings:
\begin{equation}
\label{eq:blmur}
\Exp(\rho)\approx\sqrt{\frac{\pi r n}{2(r-1)}} .
\end{equation}

More recently, Bailey et al.~\cite[Appendix~B]{cryptoeprint:2009:541} employed an uneven partitioning with probabilities $p_i$, $1\le i\le r$, for the Pollard rho method. Again under the assumption about the randomness of the restrictions of the iterating function, they provided the following heuristic result:
\begin{equation}
 \label{eq:bailey}
\Exp(\rho)\approx\sqrt{\frac{\pi n}{2(1-\sum_{i=1}^r p_i^2 )}} ,
\end{equation}
which agrees with \eqref{eq:blmur} when all $p_i$ are equal.
Combining equations~(\ref{eq:vOW}) and~(\ref{eq:bailey}), since the probability that a collision
is good is $\mathfrak{p} = 1/2$, would lead to a conjectured expected value of $\alpha$
of $ \sqrt{ \pi n / (1 - \sum_{i=1}^r p_i^2 ) }$. Theorem~\ref{th:mean} proves this result.

%%%%%%%%%%%%%%%%%%%%%%%%%%%%%%%%%%%%%%%%%%%%%%%%%%%%%%%%%%%%%%%%%%%%%%%%%%%%%%%
\subsection{Issues Caused by Uneven Partitioning}
When some partitions are used more often than others, walks become less likely to collide. Indeed, a collision involves two edges coming from two different partitions into the same node. Since every node has exactly one outgoing edge, uneven partitioning implies uneven distribution of edges among their types, and hence it becomes less likely to pick two edges of different types. This aspect is studied in the theoretical analysis below.

Another issue caused by uneven partitioning is that walks lose their mixing property, namely they behave less like random mappings than with even partitioning. This aspect is not accounted by our theoretical model, but it is discussed in Section~\ref{sec:mixing}.

%%%%%%%%%%%%%%%%%%%%%%%%%%%%%%%%%%%%%%%%%%%%%%%%%%%%%%%%%%%%%%%%%%%%%%%%%%%%%%%
\subsection{Theoretical Model of the Algorithm}
We now define an algorithm $\mathcal{A}_{\pi}$ that closely resembles $\mathcal{A}$. The only differences between $\mathcal{A}_{\pi}$ and $\mathcal{A}$ are that the walk is implemented using random permutations, and that there is no loop detection (to simplify the proof in the next section we assume that walks never loop  before they hit a distinguished node). 
Walks for $\mathcal{A}_{\pi}$ are defined as follows. Let $h_1$, \dots, $h_r$ be random permutations on $X$ such that $h_i(z)\ne z$ and $h_i(z)\ne h_j(z)$ for all $z\in X$ and $i\ne j$. 
Walks are now defined using the map
\begin{align*}
\psi_\pi\colon X &\to X\\
z &\mapsto h_{v(z)}(z) .
\end{align*}
Algorithm $\mathcal{A}_{\pi}$ is obtained from $\mathcal{A}$ by replacing line~\ref{state:l4} of the client Algorithm~\ref{alg:Aclient} with $z\leftarrow h_{i}(z)$ and deleting lines~\ref{state:l7}--\ref{state:l10}. Because of the nature of the walks, Algorithm~$\mathcal{A}_{\pi}$ does not solve the GAIP.

%%%%%%%%%%%%%%%%%%%%%%%%%%%%%%%%%%%%%%%%%%%%%%%%%%%%%%%%%%%%%%%%%%%%%%%%%%%%%%%
\subsection{Running Time of the Theoretical Model}
\label{sec:randperm}

We now state the expected running time of Algorithm~$\mathcal{A}_{\pi}$. This is essentially the same result as given in Appendix~B of Bailey et al.~\cite{cryptoeprint:2009:541}, although their work is for the Pollard rho discrete logarithm problem, whereas we are considering a slightly different situation. We also give a Heuristic~\ref{he:var}, for the standard deviation of the running time.

\begin{theorem}
\label{th:mean}
Let $n$ be the cardinality of the set $X$, $\theta$ the probability of a node being distinguished and $p_1$, \dots, $p_r$ the probabilities of choosing among $r$ random permutations on $X$. Then the number $\alpha_\pi$ of nodes visited, with repetition, before Algorithm~$\mathcal{A}_{\pi}$ terminates, has the following expected value:
\[
\Exp(\alpha_\pi)=\sqrt{\frac{\pi n}{d}}+\frac{2}{\theta}+O(\ln^4 (n)) ,
\]
where $d$ is the expected in-degree of a visited node excluding the edge used to arrive at this node\footnote{The term in-degree refers to a graph with the set of vertices $X$ and the edges $(z,\psi_\pi(z))$. For a visited vertex, the number of used incoming edges equals zero if it is a randomized starting vertex, or one otherwise.}: 
\begin{equation}
\label{eq:d}
d=1-(1-\theta)\sum_{i=1}^r p_i^2 . 
\end{equation}
\end{theorem}

\begin{proof}
We sketch an outline of the proof and refer to Stolbunov~\cite{StolbunovThesis} for the details.
The proof uses the approach of Blackburn and Murphy~\cite{BlaMur98}.  The main task is to determine the expected number of elements sampled before the first good collision. It is then standard that $1/ \theta$ further steps are required to detect a collision. Note that two collisions are expected in total.

Let $\Lambda\subset X$ denote the set of elements already visited at some stage during the execution of Algorithm~$\mathcal{A}_\pi$. For each element $z \in \Lambda$ (except for the starting point) let $z_0 \in \Lambda$ be the previous element in the walk, and suppose $z_0$ lies in partition $i$, so that $z = h_i( z_0 )$.
Let $j \in \{ 1, \dots, r \} \setminus \{ i \}$. There is an incoming edge to $z$ corresponding to partition $j$ if and only if
$h_j^{-1} ( z )$ lies in partition $j$.  Under the assumption that the partitions are random, this occurs with probability $p_j$.
Hence, the expected number of edges into $z$ coming from partition $j$ is $p_j$.
Now, since all the permutations are random and independent, the expected number of incoming edges to $z$ is the sum of the expectations for each individual permutation:
\[ \sum_{\substack{j=1\\ j \ne i}}^r p_j . \]
Now, summing over all possible choices for $i$ (given that each arises with probability $p_i$) gives
\[\sum_{i=1}^r p_i\sum_{\substack{j=1\\j\ne i}}^r p_j=
   \sum_{\substack{1\le i,j\le r\\i\ne j}}p_i p_j=
    1-\sum_{i=1}^r p_i^2 .
\]
This is the expected number of external incoming hops, for a random non-initial element of $\Lambda$. Since the proportion of initial elements equals $\theta$, hence equation~\eqref{eq:d}. 

The expected number of elements sampled to get a collision is $\sqrt{ \pi n / (2d)}$ by the 
same arguments as used by Brent-Pollard and Blackburn-Murphy.
However, a collision is only a good collision with probability $1/2$ so, using the logic behind equation~\eqref{eq:vOW}, one gets  the formula $\sqrt{ \pi n / d}$. \hfill $\Box$
\end{proof}

Note that the value $d$ in Theorem~\ref{th:mean} can easily be computed for small $r$ and known $p_i$.
When all $p_i = 1/r$ and $\theta$ tends to zero, then $d$ tends to
\[1-r\frac{1}{r^2}=\frac{r-1}{r} .\]
Hence,  Theorem~\ref{th:mean} agrees with previous results on the Pollard rho algorithm when
using $r$ partitions all of the same size, cf. \eqref{eq:blmur}. 

\begin{heuristic}
\label{he:var}
Let $\alpha_\pi$, $n$, $\theta$ and $d$ be as in Theorem~\ref{th:mean}. Then the variance of the random variable $\alpha_\pi$ approximates as:
\begin{equation}
\label{eq:varalpha}
\Var(\alpha_\pi)\approx\frac{(4-\pi)n}{d}+\frac{4-2\theta}{\theta^2}+\frac{1}{\theta}\sqrt{\frac{\pi n}{d}} .
\end{equation}
\end{heuristic}
We provide a brief argument for Heuristic~\ref{he:var} below and refer to Stolbunov~\cite{StolbunovThesis} for the details.

The total number of visited nodes $\alpha_\pi$ is the sum of the number of unique visited nodes $\lambda_\pi$ and the number $\delta_\pi$ of nodes visited twice or more. Hence 
\[\Var(\alpha_\pi) = \Var(\lambda_\pi) + \Var(\delta_\pi)+ 2 \Cov(\lambda_\pi,\delta_\pi) ,\]
where the summands correspond to the ones in~\eqref{eq:varalpha}. The probability distribution of $\lambda_\pi$ can be approximated by the (continuous) Rayleigh distribution~\cite{Rayleigh} with the following probability density function and variance:
\[ f_{\lambda_\pi}(x) \approx \frac{xd}{2n} e^{-\frac{x^2 d}{4n}} , \qquad  \Var({\lambda_\pi}) \approx \frac{(4-\pi)n}{d} . \]
When it comes to the duplicate visited nodes, chasing the good-collision distinguished node can be described as a sequence of Bernoulli trials with success probability $\theta/2$, because only half of the collisions are good. The number of trials~$\delta_\pi$ needed to get one success conforms to the geometric distribution \cite[\S6.1.2]{MR2044140}. Hence the probability mass function and the variance of $\delta_\pi$ are
\[ f_{\delta_\pi}(x)=\frac{\theta}{2}\left(1-\frac{\theta}{2}\right)^{x-1} , \qquad \Var(\delta_\pi) = \frac{4-2\theta}{\theta^2} . \]
The covariance of $\lambda_\pi$ and $\delta_\pi$ is computed using the formula (see~\cite{StolbunovThesis})
\[ \Cov(\lambda_\pi,\delta_\pi) = \Exp(\lambda_\pi\delta_\pi)-\Exp(\lambda_\pi)\Exp(\delta_\pi). \]

%%%%%%%%%%%%%%%%%%%%%%%%%%%%%%%%%%%%%%%%%%%%%%%%%%%%%%%%%%%%%%%%%%%%%%%%%%%%%%%
\subsection{Running Time Calculations}
Let the partitioning probabilities $p_1$, \dots, $p_r$ be chosen from a geometric progression with common ratio $w$ (cf. Section~\ref{sec:better-choices}). Table~\ref{tab:meanvartheory} lists the values $d$, the expected values and the standard deviations of $L_{\pi}$ for $n=2^{80}$ and $\theta=n^{-1/4}$. Mantissas are rounded to four decimal digits. 

\begin{table}
\begin{center}
\scriptsize\begin{tabular}{llllllllll}
\toprule
& \multicolumn{3}{c}{$w=1$} & \multicolumn{3}{c}{$w=1/2$} & \multicolumn{3}{c}{$w=1/4$}\\
\cmidrule(lr){2-4}\cmidrule(lr){5-7}\cmidrule(lr){8-10}
$r$&$d$&$\Exp(L_{\pi})$&$\Stdev$&$d$&$\Exp(L_{\pi})$&$\Stdev$&$d$&$\Exp(L_{\pi})$&$\Stdev$\\
\midrule
$3$ & $0.6667$ & $2.1708$ & $1.1347$ & $0.5714$ & $2.3447$ & $1.2256$ & $0.3810$ & $2.8717$ & $1.5011$\\
$4$ & $0.7500$ & $2.0467$ & $1.0698$ & $0.6222$ & $2.2470$ & $1.1746$ & $0.3953$ & $2.8191$ & $1.4736$\\
$5$ & $0.8000$ & $1.9817$ & $1.0359$ & $0.6452$ & $2.2067$ & $1.1535$ & $0.3988$ & $2.8066$ & $1.4671$\\
$6$ & $0.8333$ & $1.9416$ & $1.0149$ & $0.6561$ & $2.1882$ & $1.1438$ & $0.3997$ & $2.8035$ & $1.4655$\\
%$7$ & $0.8571$ & $1.9145$ & $1.0007$ & $0.6614$ & $2.1794$ & $1.1392$ & $0.3999$ & $2.8028$ & $1.4651$\\
%$8$ & $0.8750$ & $1.8948$ & $0.9905$ & $0.6641$ & $2.1751$ & $1.1370$ & $0.4000$ & $2.8026$ & $1.4650$\\
%$9$ & $0.8889$ & $1.8800$ & $0.9827$ & $0.6654$ & $2.1729$ & $1.1358$ & $0.4000$ & $2.8025$ & $1.4649$\\
$10$ & $0.9000$ & $1.8683$ & $0.9766$ & $0.6660$ & $2.1719$ & $1.1353$ & $0.4000$ & $2.8025$ & $1.4649$\\
%$11$ & $0.9091$ & $1.8590$ & $0.9717$ & $0.6663$ & $2.1713$ & $1.1350$ & $0.4000$ & $2.8025$ & $1.4649$\\
%$12$ & $0.9167$ & $1.8513$ & $0.9677$ & $0.6665$ & $2.1711$ & $1.1349$ & $0.4000$ & $2.8025$ & $1.4649$\\
%$13$ & $0.9231$ & $1.8448$ & $0.9643$ & $0.6666$ & $2.1709$ & $1.1348$ & $0.4000$ & $2.8025$ & $1.4649$\\
%$14$ & $0.9286$ & $1.8394$ & $0.9615$ & $0.6666$ & $2.1709$ & $1.1348$ & $0.4000$ & $2.8025$ & $1.4649$\\
%$15$ & $0.9333$ & $1.8347$ & $0.9590$ & $0.6666$ & $2.1708$ & $1.1347$ & $0.4000$ & $2.8025$ & $1.4649$\\
$16$ & $0.9375$ & $1.8306$ & $0.9569$ & $0.6667$ & $2.1708$ & $1.1347$ & $0.4000$ & $2.8025$ & $1.4649$\\
\bottomrule
\end{tabular}
\end{center}
\caption{The values $d$, $\Exp(L_{\pi})$ and $\Stdev(L_{\pi})$, when $r$ partitions are used and partitioning probabilities decrease with ratio $w$. $n=2^{80}$ and $\theta=2^{-20}$.}
\label{tab:meanvartheory}
\end{table}

The values of $d$ in the first column of Table~\ref{tab:meanvartheory} agree with $(r-1)/r$ as expected. Note also that the values of $E( L_\pi )$ in the first column converge to the expected asymptotic value of $\sqrt{\pi}  \approx 1.7724$. The values in the $w = 1/4$ column do not change significantly when $r$ is large; this is because the higher primes are used with such extremely low probability that they have no effect on the algorithm. The values in Table~\ref{tab:meanvartheory} will be used later to give an estimate of  the running time of our improved variant of the algorithm.

%%%%%%%%%%%%%%%%%%%%%%%%%%%%%%%%%%%%%%%%%%%%%%%%%%%%%%%%%%%%%%%%%%%%%%%%%%%%%%%
%%%%%%%%%%%%%%%%%%%%%%%%%%%%%%%%%%%%%%%%%%%%%%%%%%%%%%%%%%%%%%%%%%%%%%%%%%%%%%%
\section{Comparing Theory and Practice}
\label{sec:prac}

There are many reasons why we do not expect the practical Algorithm~$\mathcal{A}$ to behave as well as the theoretical Algorithm~$\mathcal{A}_{\pi}$. The aim of this section is to briefly mention one of these issues, and to develop a plausible set of heuristics for the running time of Algorithm~$\mathcal{A}$.

%%%%%%%%%%%%%%%%%%%%%%%%%%%%%%%%%%%%%%%%%%%%%%%%%%%%%%%%%%%%%%%%%%%%%%%%%%%%%%%
\subsection{Mixing of Adding Walks}
\label{sec:mixing}

As is standard, the theoretical analysis assumes truly random walks.  However, we are using adding walks in a group, and such walks are not close to uniformly distributed if they are short. The mixing time is a measure of how long a walk runs before its values start to appear uniformly distributed. It is beyond the scope of this paper to analyse such issues in detail. We mention that Dai and Hildebrand~\cite{MR1483024} have studied the mixing time of adding walks. They show that adding walks on $r$ partitions need a slack of $O(n^{(2/(r-1))+\epsilon})$ hops before they converge to the uniform distribution. 

However, it is worth noting that Algorithm~$\mathcal{A}$ does not necessarily need walks to be uniformly distributed after a certain number of hops. Instead it needs walks to collide. Just because walks have not yet reached uniform sampling does not prevent collisions from occurring.

%%%%%%%%%%%%%%%%%%%%%%%%%%%%%%%%%%%%%%%%%%%%%%%%%%%%%%%%%%%%%%%%%%%%%%%%%%%%%%%
\subsection{Experiments}
\label{sec:practical}

To get a better idea of how the algorithm works in practice, we have performed a suite of experiments. We report one of them in this paper and refer to Stolbunov~\cite{StolbunovThesis} for more details.

Our numerical experiments are for $X = G$ (i.e., $G$ acting on itself) being an abstract group of the form $\mathbb{Z}_{n_1}\oplus\dots\oplus \mathbb{Z}_{n_s}$, where $n_{i+1}\mid n_i$ and $n_i\ge2$ for all~$i$. The integer $s$ is the \emph{rank} of $G$. The supporting set is randomly chosen, though it is checked that it generates the group.

For calculations we use a Linux cluster of $32$ quad-core Intel X5550 processors clocked at 2.67~GHz. The code is written in C++. We use a single-threaded implementation of Algorithm~$\mathcal{A}$, such that one thread alternates between $x_0$- and $x_1$-walks. The same experiment is run on all CPU cores in parallel but with different random generator seeds. 

Group elements are represented by arrays of 64-bit integers. We make use of a hash function $\Hash : G\to \{0,1\}^{32}$ implemented using the 64 to 32 bit hash function of Wang~\cite{Wang}. The partitioning function $v(z)$ is computed by reducing $\Hash(z)$ modulo a sufficiently large integer whose residues can be partitioned with the correct proportions. Wang's hash function uses bit shifts, negations, additions and XOR operations. This helps to make sure that $v(z)$ and $v(\psi(z))$ look like independent random variables, which is important because correlations between the functions $\psi(z)$ and $v(z)$ can result in undesirable loops in the walk.

Let $\theta$ be the desired distinguished  point probability.
We declare an element $z$ to be distinguished iff $\Hash(z)\equiv 0 \mod \lfloor 1/\theta \rceil$, where $\lfloor \cdot \rceil$ is the rounding to the nearest integer. Although Algorithm~$\mathcal{A}$ has polynomial memory requirements, we find it practical to use an $O(n^{1/4})$ amount of storage\footnote{Let us justify the suitability of this choice by an example. Suppose one tries to solve a $\Cl$-GAIP over a $244$-bit field, a problem size proposed for isogeny-based cryptosystems~\cite{MR2654134}. Since the group size (i.e., class number) $n\approx2^{122}$, the database of distinguished nodes should store $L\theta \sqrt{n}$ nodes, which is less than $2^{33}$ on average. Since the class number is approximately $122$ bits long, one entry of the database (binary tree) of distinguished nodes would occupy $48$ bytes, of which $16$ bytes are used by a hashed $j$-invariant, $16$ bytes by a compressed class group element and $16$ bytes by two pointers. The whole database would occupy not more than $384$~gigabytes of disk space, which we find to be quite moderate.}, namely to choose
\[
\theta=n^{-\frac{1}{4}} .
\]
This is compatible with the work of Schulte-Geers~\cite{schulte2000}. The database of distinguished nodes is implemented as a binary tree.

For the starting randomization of walks we use the 64-bit Mersenne twister pseudorandom generator~\cite{MatsumotoN98}. A pseudorandom element $g_r\in G$ acts on the initial node to create the starting point of the new walk.

%%%%%%%%%%%%%%%%%%%%%%%%%%%%%%%%%%%%%%%%%%%%%%%%%%%%%%%%%%%%%%%%%%%%%%%%%%%%%%%
\subsection{Choosing the Number of Experiments}
\label{sec:numexp}
Let $k$ be the number of experiments and $L_k$ the average value of $L$ over $k$ experiments. According to the central limit theorem~\cite[\S7.2.1]{MR2044140}, the probability distribution of the random variable $L_k$ approaches the normal distribution with the mean $\Exp(L)$ and the variance $\Var(L)/k$ as $k$ approaches infinity. For the normal distribution,  over $99.7~\%$ of the values lie within three standard deviations away from the mean. Thus, assuming $k$ is big enough, we have that
\[
\Pr\left[ L_k-3\frac{\Stdev(L)}{\sqrt{k}} \le E(L) \le L_k+3\frac{\Stdev(L)}{\sqrt{k}} \right] > 0.997 .
\]
When measuring $\Exp(L)$, we use two levels of accuracy: the result lies within $\pm0.1~\%$ of the true value for the experiments satisfying $\log (n)\le 44$, and within $\pm0.5~\%$ of the true value otherwise. Thus we can use the inequalities
\begin{equation}
\label{eq:kformulae}
k_{1}\ge\left(\frac{3\Stdev(L)}{0.001\Exp(L)}\right)^2, \qquad k_{2}\ge\left(\frac{3\Stdev(L)}{0.005\Exp(L)}\right)^2
\end{equation}
to find the sufficient number of experiments for the two accuracy levels. For a preliminary estimation of the number of experiments we use the formulae for $\Exp(L)$ and $\Stdev(L)$ obtained in Section~\ref{sec:randperm}. This gives us the values
\[
k_{1}=2459137, \qquad  k_{2}=98368\ ,
\]
computed as maximums over all possible parameters in Experiment~\ref{exp:1}.

Our experiments have shown that, in most cases, both the sample mean and the sample standard deviation differ from the results of Theorem~\ref{th:mean} by approximately the same factor, which cancels out in~\eqref{eq:kformulae}. This means that the obtained numbers $k_1$ and $k_2$ fit for the probability distributions under observation.

%%%%%%%%%%%%%%%%%%%%%%%%%%%%%%%%%%%%%%%%%%%%%%%%%%%%%%%%%%%%%%%%%%%%%%%%%%%%%%%
\subsection{Experimental Measurement of $L$}
In this section we measure $\Exp(L)$ by means of experimentation and assemble results in a table so that they can be used for arbitrary GAIP instances in the future.

\begin{experiment}[Measuring $L$ in Arbitrary Groups]
\label{exp:1}
For each of the values\footnote{We use $n>2^{27}$ because otherwise $L$ is highly affected by looped walks: every loop increases the number of visited nodes by $30 n^{1/4}$.} $\lceil\log (n)\rceil\in\{28, 32, 36, \dots, 56\}$, $r\in\{3,4, \dots, 16\}$ and $w\in\{1$, $3/4$, $1/2$, $1/3$, $1/4\}$ conduct a set of $k_1$ ($k_2$ for $n>2^{44}$) experiments. In each experiment choose a random\footnote{For each $m\in\{28, 32, 36, \dots, 56\}$ we sample uniformly from the set of isomorphism classes of abelian groups of order $n$ and rank at most $r$, where $2^{m-1} + 1 \le n \le 2^m$.} group $G$ and a random subset of $r$ elements that generates~$G$. Use $\theta=n^{-1/4}$ and the partitioning probabilities decreasing with ratio $w$. 
\end{experiment}

\begin{table}
\begin{center}
\footnotesize\begin{tabular}{llllllllll}
\toprule
&& \multicolumn{8}{c}{$\lceil \log (n) \rceil$} \\
\cmidrule{3-10}
$w$ & $r$ & $28$ & $32$ & $36$ & $40$ & $44$ & $48$ & $52$ & $56$\\
\midrule
\multirow{6}{*}{$1$}
 & $3$ & $2.8547$ & $2.9982$ & $3.1380$ & $3.2735$ & $3.4079$ & $3.5355$ & $3.6681$ & $3.7812$ \\
 & $4$ & $2.2923$ & $2.3101$ & $2.3247$ & $2.3371$ & $2.3484$ & $2.3518$ & $2.3661$ & $2.3661$ \\
 & $5$ & $2.1039$ & $2.0975$ & $2.0968$ & $2.0984$ & $2.0978$ & $2.1007$ & $2.1009$ & $2.1004$ \\
 & $6$ & $2.0178$ & $2.0099$ & $2.0052$ & $2.0032$ & $2.0026$ & $2.0038$ & $2.0022$ & $2.0023$ \\
 & $10$ & $1.9021$ & $1.8932$ & $1.8862$ & $1.8849$ & $1.8831$ & $1.8816$ & $1.8769$ & $1.8879$ \\
 & $16$ & $1.8575$ & $1.8455$ & $1.8407$ & $1.8384$ & $1.8361$ & $1.8302$ & $1.8369$ & $1.8357$ \\
\midrule
\multirow{6}{*}{$\frac{1}{2}$}
 & $3$ & $3.1089$ & $3.2761$ & $3.4406$ & $3.5985$ & $3.7500$ & $3.9086$ & $4.0331$ & $4.1785$ \\
 & $4$ & $2.6071$ & $2.6436$ & $2.6723$ & $2.6938$ & $2.7101$ & $2.7307$ & $2.7315$ & $2.7406$ \\
 & $5$ & $2.4586$ & $2.4665$ & $2.4723$ & $2.4782$ & $2.4802$ & $2.4875$ & $2.4821$ & $2.4776$ \\
 & $6$ & $2.4000$ & $2.4022$ & $2.4063$ & $2.4068$ & $2.4086$ & $2.4079$ & $2.4069$ & $2.4128$ \\
 & $10$ & $2.3529$ & $2.3536$ & $2.3527$ & $2.3553$ & $2.3563$ & $2.3486$ & $2.3533$ & $2.3557$ \\
 & $16$ & $2.3516$ & $2.3523$ & $2.3519$ & $2.3519$ & $2.3536$ & $2.3524$ & $2.3465$ & $2.3576$ \\
\midrule
\multirow{6}{*}{$\frac{1}{4}$}
 & $3$ & $3.8596$ & $4.1194$ & $4.3652$ & $4.5978$ & $4.8213$ & $5.0395$ & $5.2484$ & $5.4338$ \\
 & $4$ & $3.5425$ & $3.6694$ & $3.7582$ & $3.8280$ & $3.8771$ & $3.9015$ & $3.9372$ & $3.9517$ \\
 & $5$ & $3.4753$ & $3.5732$ & $3.6423$ & $3.6830$ & $3.7103$ & $3.7322$ & $3.7295$ & $3.7407$ \\
 & $6$ & $3.4608$ & $3.5566$ & $3.6145$ & $3.6526$ & $3.6743$ & $3.6985$ & $3.6845$ & $3.6845$ \\
 & $10$ & $3.4578$ & $3.5486$ & $3.6064$ & $3.6454$ & $3.6658$ & $3.6672$ & $3.6853$ & $3.6833$ \\
 & $16$ & $3.4607$ & $3.5498$ & $3.6070$ & $3.6427$ & $3.6639$ & $3.6747$ & $3.6808$ & $3.6880$ \\
\bottomrule
\end{tabular}
\end{center}
\caption{Expected values of $L$ obtained experimentally for certain choices of $r$ and $w$.}
\label{tab:Larb}
\end{table}

A subset of results is listed in Table~\ref{tab:Larb}, where mantissas are rounded to four decimal digits. Full data for $3 \le r \le 16$ and $w \in \{1, 3/4, 1/2 , 1/3, 1/4 \}$ are available in \cite{StolbunovThesis}. The entire experiment took 51 days of parallel processing on 128 cores.

When $w = 1$ and $r = 16$ one sees good agreement between Table~\ref{tab:Larb} and Table~\ref{tab:meanvartheory}, which suggests that our implementation is working well. In other cases we see that $L$ is significantly larger than $L_\pi$, which shows that the theoretical analysis is over-optimistic about the behaviour of these pseudorandom walks. The results also confirm that $r = 3$ is not a good choice in practice.

\begin{figure}
\centering
\includegraphics[width=\textwidth]{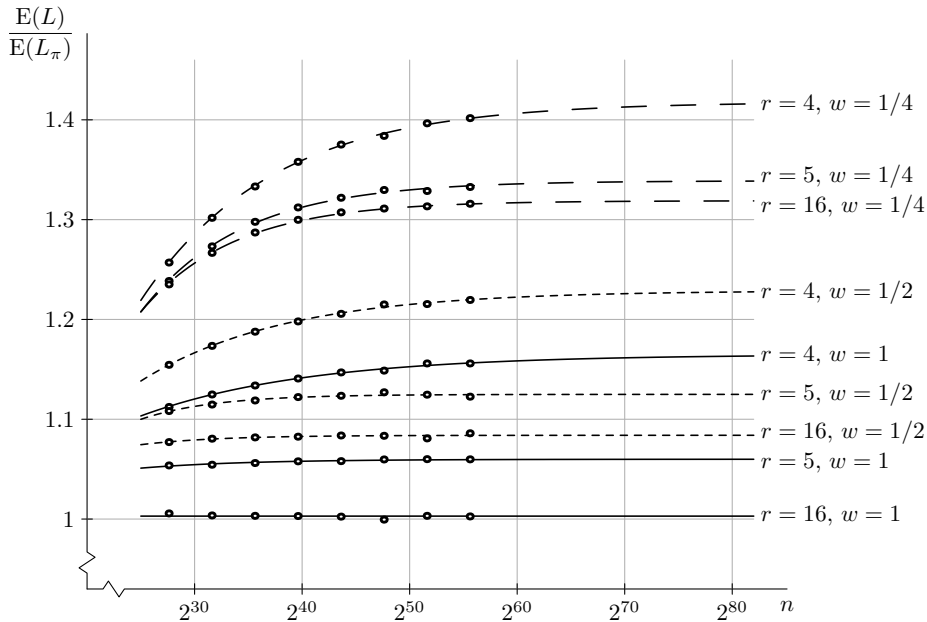}
\caption{Values of $\sigma = \Exp(L) / \Exp(L_{\pi}) $ obtained experimentally and their approximations extended to $n=2^{80}$.}
\label{fig:differences}
\end{figure}

Figure~\ref{fig:differences} graphs some values of the practice-to-theory ratio 
\[\sigma = \frac{\Exp(L)}{\Exp(L_{\pi})}.\]
Round dots depict our experimental results, and lines are their approximating functions (solid lines are $w = 1$, short-dashed lines are $w = 1/2$ and long-dashed lines are $w = 1/4$). For a fixed $w$, values of $\sigma$ for $5 < r < 16$ lie between $r=5$ and $r=16$. One can observe an increased roughness of experimental results for $n>2^{44}$ due to the increased confidence interval. The graphs suggest that, for $r > 3$, the difference between $\Exp(L)$ and $\Exp( L_{\pi} )$ is fairly stable as $n$ grows. Hence, when $r > 3$ we feel confident extrapolating actual values for $\Exp(L)$ from our formulae for $\Exp( L_\pi )$ and the experimentally determined correction factors $\sigma$.

\begin{remark}
Recently Montenegro proposed a heuristic for estimating the number of hops in birthday attacks~\cite{MontenegroHeuristic}. His idea is to estimate the probability of short cycles, i.e. if two walks (with independent partitioning functions) are started from the same position, then what is the probability that they intersect soon? The lower this probability is, the sooner the algorithm will terminate. Applied to adding walks in an abelian group, this means that if two walks include short subsequences of edges which are equivalent up to the order of edges, these subsequences do not change the relative position of these walks. Although Montenegro only gives examples for Pollard's and Teske's walks, his heuristic also applies to walks with uneven partitioning. The probability $P_1$ that two independent walks started from $x_0=y_0$ have a collision after one hop equals
\[P_1=\Pr\left[x_1=y_1\right]=\sum_{i=1}^r p_i^2.\]
If we only consider collisions after one hop, then Montenegro's heuristic gives an approximation similar to what we obtained in Theorem~\ref{th:mean}:
\[\Exp(\lambda)\approx\sqrt{\frac{\pi n}{1-P_1}}.\]
The probability $P_2$ that a collision occurs on the second hop is
\[P_2=\Pr\left[(x_1\ne y_1) \land (x_2=y_2)\right]=(1-P_1)P_1^2,\]
and Montenegro's heuristic gives
\begin{equation}
\label{eq:montenegro}
\Exp(\lambda)\approx\sqrt{\frac{\pi n}{1-(P_1+P_2)}}=\sqrt{\frac{\pi n}{1-P_1-P_1^2+P_1^3}}.
\end{equation}
The calculation can be continued to more hops, but since probabilities of collisions become much smaller than $P_2$, this will result in very small numerical changes. 

We have calculated the expected values of $L$ using \eqref{eq:montenegro} and found that for $r\ge 6$ the heuristic agrees pretty well with our practical results, giving only up to $3.4~\%$ error for $w=1/2$ and up to $5.6~\%$ error for $w=1/4$.
\end{remark}

%%%%%%%%%%%%%%%%%%%%%%%%%%%%%%%%%%%%%%%%%%%%%%%%%%%%%%%%%%%%%%%%%%%%%%%%%%%%%%%
%%%%%%%%%%%%%%%%%%%%%%%%%%%%%%%%%%%%%%%%%%%%%%%%%%%%%%%%%%%%%%%%%%%%%%%%%%%%%%%
\section{The Algorithm in Practice}
\label{sec:predictions}

We now discuss how the isogeny algorithm performs in practice. We focus on the case of ideal class groups of maximal orders in CM fields coming from $\End(E)$ where $E$ is a randomly chosen elliptic curve over $\F_p$ and $p$ is a randomly chosen $160$-bit prime. 
We also speculate on how the algorithm will perform for larger fields at the end of this section.

We have already obtained a good theoretical and experimental understanding of the algorithm for the group action problem. It is necessary now to include the cost of computing isogenies. The next section gives some estimates for the running time of computing isogenies of prime degree.

%%%%%%%%%%%%%%%%%%%%%%%%%%%%%%%%%%%%%%%%%%%%%%%%%%%%%%%%%%%%%%%%%%%%%%%%%%%%%%%
\subsection{Cost of Computing Isogenies}
\label{sec:cost_isogenies}

Consider the cost of computing the action by a prime ideal in the isogeny graph.
One has an elliptic curve and an ideal of norm $\ell$.
One must factor the modular polynomial to determine the possible $j$-invariants of $\ell$-isogenous curves,
one must perform Elkies' algorithm to determine the kernel polynomials for these isogenies, and then
one must use the technique from~\cite{MR1975526} to determine which is the correct kernel and hence
which is the correct isogeny\footnote{If two or more consecutive hops are made by the same split isogeny degree $\ell$, and there are no vertical $\ell$-isogenies, then it is sufficient to choose the correct isogeny only at the first hop. On each subsequent hop one simply checks that the $j$-invariant does not match the previous one. This provides extra saving, especially when the partitioning is uneven. This extra saving is not accounted in Table~\ref{tab:total_timings}.}.  It is not necessary to apply V{\' e}lu's formulae at this stage.
We assume the modular polynomials have been precomputed and reduced to the finite field~$\F_q$.
Since the modular polynomial has $O( \ell^2 )$ coefficients one performs  
%$O( \ell^2  \log(q)^2  )$ bit operations 
$O( \ell^2 )$ field operations
to evaluate the modular polynomial at the target $j$-invariant.
An expected 
%$O( \ell^2 \log( \ell) \log( q )  )$ 
$O( \ell \log( \ell ) \log( q )  )$ field operations are performed to find the roots of the polynomial, employing fast polynomial arithmetic. Finally, $O( \ell^2 )$ field operations are used by Elkies' algorithm.
Hence one expects the time of one $\ell$-hop to grow like 
\begin{equation}
\label{eq:hoptime}
O\left( \ell^{2} + \ell \log( \ell ) \log( q ) \right)
\end{equation}
field operations.

\begin{figure}
\hfill
\scriptsize\begin{tabular}{ll}
\toprule
$\ell$ & Time, s \\
\midrule
$3$ & $0.002870$ \\
$5$ & $0.004799$ \\
$7$ & $0.006898$ \\
$11$ & $0.012113$ \\
$13$ & $0.015261$ \\
$17$ & $0.022376$ \\
$19$ & $0.026499$ \\
$23$ & $0.036209$ \\
$29$ & $0.052346$ \\
$31$ & $0.058230$ \\
$37$ & $0.084434$ \\
$41$ & $0.092049$ \\
$43$ & $0.106742$ \\
$47$ & $0.116152$ \\
$53$ & $0.143732$ \\
$59$ & $0.150925$ \\
\bottomrule
\end{tabular}
\hfill
\scriptsize\begin{tabular}{ll}
\toprule
$\ell$ & Time, s \\
\midrule
$61$ & $0.168179$ \\
$67$ & $0.225253$ \\
$71$ & $0.225503$ \\
$73$ & $0.249365$ \\
$79$ & $0.270537$ \\
$83$ & $0.284242$ \\
$89$ & $0.305786$ \\
$97$ & $0.341988$ \\
$101$ & $0.353268$ \\
$103$ & $0.362195$ \\
$107$ & $0.375111$ \\
$109$ & $0.384550$ \\
$113$ & $0.403007$ \\
$127$ & $0.467993$ \\
$131$ & $0.579427$ \\
$137$ & $0.624039$ \\
\bottomrule
\end{tabular}
\hfill
\begin{minipage}[c]{.5\textwidth}
\includegraphics[width=\textwidth]{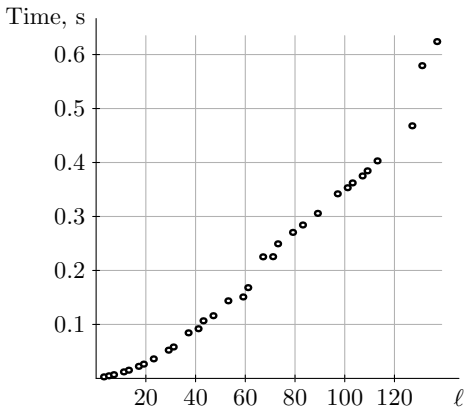}
\end{minipage}
\hfill
\caption{Average running time of one $\ell$-isogeny (i.e., action by a prime ideal of norm $\ell$) for elliptic curves over $160$-bit prime fields.}
\label{fig:prime_action_timing}
\end{figure}

We computed average timings using the ClassEll package by Stolbunov~\cite{classell}. The package implements the ideal class group action on sets of ordinary elliptic curves. The experiment was run on Intel X5550 processors clocked at 2.67~GHz, the code executed at approximately $6799$ millions instructions per second (MIPS). 
The data was gathered by repeatedly ($20000$ times) generating a random $160$-bit prime $p$ and a random ordinary elliptic curve over $\F_p$ with a fundamental Frobenius discriminant. The time spent on one action by a prime ideal, for prime ideals of all split norms less than or equal to $137$, was recorded. To increase the accuracy, we performed more hops for smaller primes. Results are given in Fig.~\ref{fig:prime_action_timing}. We can observe bumps when $\ell$ moves over degrees of two which is typical for the polynomial multiplication by number-theoretic transform.

%%%%%%%%%%%%%%%%%%%%%%%%%%%%%%%%%%%%%%%%%%%%%%%%%%%%%%%%%%%%%%%%%%%%%%%%%%%%%%%
\subsection{Ideal Class Groups}

In Experiment~\ref{exp:1} we used the uniform distribution of finite abelian groups. However, the structure of ideal class groups is not that random; the following observations are known as Cohen-Lenstra heuristics~\cite{MR756082}: the odd part of the class group of an imaginary quadratic field is quite rarely non-cyclic; if $p$ is a small odd prime, the proportion of imaginary quadratic fields whose class number is divisible by $p$ is close to $1/p+1/p^2$. The distribution of group structures in the isogeny problem is further affected by the fact that the imaginary quadratic orders are chosen as endomorphism rings of random elliptic curves. Nevertheless, our experiments show that the difference between values $\Exp(L)$ for random isogeny problem instances\footnote{Parameters: $\lceil\log (p)\rceil=90$, $4\le r\le 16$; $w$, $\theta$, $c_{\max}$ and $k_1$ are as in Experiment~\ref{exp:1}.} and for random GAIP instances lies within the margin of error $0.2~\%$. The same holds for the standard deviation of $L$.

Due to the numerical results of Jacobson, Ramachandran and Williams~\cite{MR2282917} we know that the average maximum norm of the prime ideals required to generate the class group of $\mathbb{Q}(\sqrt{\Delta})$ for $-10^{11}<\Delta<0$ approximately equals $0.60191\ln\lvert\Delta\rvert$, and the number of prime ideals required to generate these class groups averages at approximately $3.3136$. We assume that these results apply to our problem size as well. Hence for a random ideal class group of a $162$-bit discriminant, it is very likely that a generating set of four prime ideals with the maximum norm $67$ can be found. This observation is used in the next section where we model the choice of primes.

We make an assumption that walks with a supporting set that consists of ideals of small prime norm behave similar to walks when the supporting set consists of random group elements.

%%%%%%%%%%%%%%%%%%%%%%%%%%%%%%%%%%%%%%%%%%%%%%%%%%%%%%%%%%%%%%%%%%%%%%%%%%%%%%%
\subsection{Predicted Results}
In this section we estimate the time needed for solving a random instance of the isogeny problem over a $160$-bit finite field using various numbers of partitions $r$ and partitioning probabilities $\vec p$. The expected serial running time is computed using equation \eqref{eq:exptime}, which can be written as
\[ \sigma\,\Exp(L_{\pi})\,\sqrt{n}\, \vec p\, \vec t .\]
The values $\Exp(L_{\pi})$ are computed using Theorem~\ref{th:mean} and approximations for $\sigma$ are based on our experimental data 
(partially displayed on Fig.~\ref{fig:differences}). We take $n = 2^{80}$. What remains is to compute the average running time $\vec p\, \vec t$ of one hop.

For the isogeny problem, the supporting set $\WB$ should be chosen to consist of prime ideals above the smallest integer primes which split in $\OO$. If necessary, one or more prime ideals of larger norm are included in $\WB$ to ensure that $\WB$ generates $\Cl( \OO )$. 
To compute the average product $\vec p\, \vec t$ for given $r$ and $w$, we enumerate all subsets $\WB$ of $r$ primes larger or equal to $3$ with the $r-4$ smallest primes in $\WB$ being less than or equal to\footnote{Because approximately half of primes are split.} $\mathrm{prime}_{2r-7}$, and the largest prime in $\WB$ lying between $67$ and $\max(67,\mathrm{prime}_{2r+1})$. For every set $\WB$, a timing vector $\vec t$ is constructed using the data on Fig.~\ref{fig:prime_action_timing}. Hence we compute the average $\vec p\, \vec t$ over all~$\WB$.

In Table~\ref{tab:total_timings} we give estimated times for solving a random instance of the isogeny problem over a $160$-bit finite field (equivalently, the $\Cl$-GAIP problem in $\Cl( \OO )$ where $\OO = \End( E) $ for an elliptic curve over a $160$-bit finite field). The time is provided in years of serial execution on one Intel X5550 2.67~GHz CPU core. On a cluster with hundreds of thousands of cores the problem can be solved in a matter of hours.

\begin{table}
\begin{center}
\footnotesize\begin{tabular}{llllll}
\toprule
r\textbackslash w & $1$ & $3/4$ & $1/2$ & $1/3$ & $1/4$\\
\midrule
$4$ & $8708$ & $6940$ & $5429$ & $4727$ & $4690$ \\
$5$ & $6455$ & $4495$ & $2758$ & $1925$ & $1652$ \\
$6$ & $5514$ & $3396$ & $1755$ & $1130$ & $988$ \\
$7$ & $5068$ & $2827$ & $1334$ & $904$ & $858$ \\
$8$ & $4891$ & $2530$ & $1154$ & $847$ & $848$ \\
$9$ & $4930$ & $2415$ & $1093$ & $\mathbf{842}$ & $856$ \\
$10$ & $5549$ & $2548$ & $1110$ & $858$ & $870$ \\
$11$ & $6391$ & $2723$ & $1132$ & $874$ & $885$ \\
$12$ & $7409$ & $2915$ & $1157$ & $891$ & $903$ \\
$13$ & $8485$ & $3095$ & $1180$ & $906$ & $919$ \\
$14$ & $9519$ & $3255$ & $1205$ & $923$ & $932$ \\
$15$ & $10636$ & $3396$ & $1225$ & $937$ & $944$ \\
$16$ & $\mathbf{12200}$ & $3541$ & $1242$ & $949$ & $955$\\
\bottomrule
\end{tabular}
\end{center}
\caption{Expected serial time (years) needed to solve a random $\Cl$-GAIP over a $160$-bit field.}
\label{tab:total_timings}
\end{table}

We see from Table~\ref{tab:total_timings} that the best combination $r=9$ and $w=1/3$ is approximately $14$ times faster than $16$ equally-sized partitions (both timings are in bold). In fact all values within $7\le r\le 16$, $w\in\{1/3, 1/4\}$ provide good speeds.

For the rest of the section we briefly consider the question of how much faster our algorithm is than the GHS algorithm as $q \to \infty$. Both algorithms require $\tilde{O}( \sqrt{n} )$ bit operations, but it is not immediately clear that the ratio of running times is bounded as $q \to \infty$.  Let us compare $r=16$, $w=1$ with $r=9$, $w=1/3$. First we make a simplifying assumption: for any problem instance, a supporting set $H$ consisting of the $r-1$ smallest split primes and one prime close to $\ln (q)$, generates the class group. Using the prime number theorem we approximate primes in $H$ by $\ell_i\approx 2i\ln(2i)$, for $1\le i\le r-1$. We also approximate $\ell_r \approx \ln(q)$. Since $\ell_i<\log (q)$ for sufficiently large $q$, the complexity \eqref{eq:hoptime} of one $\ell$-hop is $O( \ell \ln( \ell ) \ln( q ) )$ field operations, which we further approximate by $c\, \ell \ln (q)$ for some constant $c$. The improvement ratio (i.e., speedup) is
\begin{multline*}
\frac{\left(\Exp(L)\, \vec p\, \vec t\right)\rvert_{\begin{subarray}{l}r=16\\w=1\end{subarray}}}{\left(\Exp(L)\, \vec p\, \vec t\right)\rvert_{\begin{subarray}{l}r=9\\w=1/3\end{subarray}}} \approx 
\frac{1.836}{3.023}\ 
\frac{\frac{1}{16}\sum_{i=1}^{15} 2i \ln (2i) \ln (q)+\frac{1}{16}\ln^2 (q)}
{\frac{6561}{9841}\sum_{i=1}^{8}\left(\frac{1}{3}\right)^{i-1} 2i \ln (2i) \ln (q)+\frac{1}{9841}\ln^2 (q)}  \\
\approx 0.607\ 
\frac{44.046+\frac{1}{16}\ln (q)}
{3.682+\frac{1}{9841}\ln (q)}
\to 0.607\ \frac{9841}{16}\approx 373 \qquad \text{as } q \to \infty .
\end{multline*}
Hence the improvement ratio slowly grows with $q$ and stabilizes at few hundreds for a very large $q$ (at $\ln(q)>2^{25}$ in the example above). Sure, problems of that size are not feasible, and $9$ primes are probably not sufficient to generate a class group that big. The growth of the improvement ratio is hard to predict, but we see no reasons for it to overcome $O(1)$ as~$q \to \infty$.

%%%%%%%%%%%%%%%%%%%%%%%%%%%%%%%%%%%%%%%%%%%%%%%%%%%%%%%%%%%%%%%%%%%%%%%%%%%%%%%
%%%%%%%%%%%%%%%%%%%%%%%%%%%%%%%%%%%%%%%%%%%%%%%%%%%%%%%%%%%%%%%%%%%%%%%%%%%%%%%
\section{Conclusion}
\label{sec:concluding}

In this paper we have improved the GHS algorithm for constructing isogenies between ordinary elliptic curves. Our improvement is by an $O(1)$ factor, which was estimated to be approximately $14$ for random $160$-bit elliptic curves with comparable conductors. This is a significant acceleration. Nevertheless, the asymptotic complexity of the $\F_q$-isogeny problem for curves with comparable conductors is $O (q^{1/4 + o(1)} \log^2 (q) \log( \log (q) ))$ field operations, as before.

%%%%%%%%%%%%%%%%%%%%%%%%%%%%%%%%%%%%%%%%%%%%%%%%%%%%%%%%%%%%%%%%%%%%%%%%%%%%%%%
%%%%%%%%%%%%%%%%%%%%%%%%%%%%%%%%%%%%%%%%%%%%%%%%%%%%%%%%%%%%%%%%%%%%%%%%%%%%%%%
\section*{Acknowledgements}
The paper was created through a collaboration of two authors whose names are listed alphabetically. The work was initiated during a two-month research visit of Anton Stolbunov to Steven Galbraith. Stolbunov would like to thank Department of Telematics, Norwegian University of Science and Technology, for the financial support of his research and that visit. We thank Gaetan Bisson and Edlyn Teske for their valuable comments on this paper. 
%%%%%%%%%%%%%%%%%%%%%%%%%%%%%%%%%%%%%%%%%%%%%%%%%%%%%%%%%%%%%%%%%%%%%%%%%%%%%%%
%%%%%%%%%%%%%%%%%%%%%%%%%%%%%%%%%%%%%%%%%%%%%%%%%%%%%%%%%%%%%%%%%%%%%%%%%%%%%%%
\bibliographystyle{plain}
\bibliography{cl-gaip}

\end{document}